\documentclass[conference]{IEEEtran}
\usepackage[utf8]{inputenc}

\usepackage{amsmath}            %equations
\usepackage{url}                %nice URL formatting
\usepackage{cite}               %ref condensing, e.g. [5-7]
\usepackage{xcolor}             %\textcolor
\usepackage{graphicx}           %figures
\graphicspath{{figures/}}

\title{U.S. Test System with High Spatial and Temporal Resolution for Renewable Integration Studies}

\author{\IEEEauthorblockN{Yixing Xu, Nathan Myhrvold, Dhileep Sivam, Kaspar Mueller, Daniel J. Olsen, Bainan Xia,
Daniel Livengood,\\Victoria Hunt, Benjamin Rouillé d'Orfeuil, Daniel Muldrew, Merrielle Ondreicka, Megan Bettilyon}\IEEEauthorblockA{Intellectual Ventures, Bellevue, Washington, USA\\Email: yxu@intven.com}}
\usepackage{tikz}
\usepackage{textcomp}
\usepackage{hyperref}
\usepackage{lipsum}

\newcommand{\copyrighttext}{%
    \footnotesize \textcopyright 2020 IEEE.  Personal use of this material is permitted.  Permission from IEEE must be obtained for all other uses, in any current or future media, including reprinting/republishing this material for advertising or promotional purposes, creating new collective works, for resale or redistribution to servers or lists, or reuse of any copyrighted component of this work in other works.}

\newcommand{\copyrightnotice}{%
    \begin{tikzpicture}[remember picture,overlay]
    \node[anchor=south,yshift=8pt] at (current page.south) {\fbox{\parbox{\dimexpr\textwidth-\fboxsep-\fboxrule\relax}{\copyrighttext}}};
    \end{tikzpicture}%
}
\date{November 2019}

\begin{document}
% Avoid ------ in bibliography for repeated authors [1/2]
%   (2/2 is in bibliography.bib)
\bstctlcite{IEEEexample:BSTcontrol}

\maketitle
\copyrightnotice

\begin{abstract}
Planning for power systems with high penetrations of variable renewable energy requires higher spatial and temporal granularity. However, most publicly available test systems are of insufficient fidelity for developing methods and tools for high-resolution planning. This paper presents methods to construct open-access test systems of high spatial resolution to more accurately represent infrastructure and high temporal resolution to represent dynamics of demand and variable resources.

To demonstrate, a high-resolution test system representing the United States is created using only publicly available data. This test system is validated by running it in a production cost model, with results compared against historical generation to ensure that they are representative. The resulting open source test system can support power system transition planning and aid in development of tools to answer questions around how best to reach decarbonization goals, using the most effective combinations of transmission expansion, renewable generation, and energy storage.
\end{abstract}


\let\thefootnote\relax\footnote{\hspace{-2pt}Funding provided by The Global Good Fund I, LLC (www.globalgood.com)}

\section{Introduction}\label{sec_intro}
Continental-scale electric power systems are some of the most complex human creations in existence. Their smooth operation is critical for modern life, and yet they are undergoing a transformation due to the need to reduce anthropogenic greenhouse gas emissions. This transformation will upend decades of assumptions about power system planning and operation, and is therefore critical to get right.

Effectively planning this transformation will require high-quality models of existing infrastructure. Although recreation of high-resolution spatial data is sometimes publicly available, detailed information is often restricted on national security grounds. This limitation led to creation of `synthetic' models of the U.S. power system by Texas A\&M University (`the TAMU network'), representative of real patterns of geography and network topology while disclosing no protected information \cite{birchfield2017_tpwrs}. Some work has been done to add time-series load data to this synthetic network \cite{li2018}, but to the best of the authors' knowledge, a full dataset providing granular temporal data for loads and variable renewable energy (VRE) for a large-scale network of the U.S. power system is still unavailable.

This paper describes methods used to create high-resolution test systems to better represent the current state of existing power systems. These efforts fall into three broad categories: collecting static data describing system infrastructure, adding high-resolution time-series data to represent the time-varying nature of demand and variable renewable resources, and refining the dataset based on operational simulations. Time-matched historical data are used for demand and VRE generation to capture their inherent correlations. This dataset, and associated methods, are the main contributions of this paper.
%The methods are described in more detail in Sections \ref{sec_static}-\ref{sec_tsd}.

As a demonstration, a high-resolution test system is built representing the U.S. power system in 2016, using entirely publicly available information. This test system is publicly available at \cite{zenodo_data}. This test system complements modern, publicly available test systems such as the TAMU networks \cite{tamu2019}, ReEDS \cite{nrel_reeds}, and RTS-GMLC \cite{rts_gmlc}; each of which have strengths and drawbacks that make them more appropriate for certain types of studies. The TAMU networks have high spatial resolution, but lack profiles for VRE generation; ReEDS has good data for long-term planning, but has low spatial resolution; and the RTS-GMLC network has high-resolution load and VRE data, but covers a small geographic area.

%Our contribution includes a set of methods used to create high-resolution test systems to better represent the current state of existing power  systems, and an open source U.S. test system with synchronized hourly time-series data.

%Although time-series demand data were added to the TAMU network in its most recent release \cite{tamu2019} (September 2019), the methods described in this paper produce time-series data for demand, wind, solar, and hydro resources. The methods described in this paper are applicable to any other network where topology is present without time-series data, or when time-series data representing different years are desired.

\section{Static Data: Network and Generator Capacities}\label{sec_static}
\subsection{Network Topology}

\begin{figure}
    \centering
    \includegraphics[width=\linewidth]{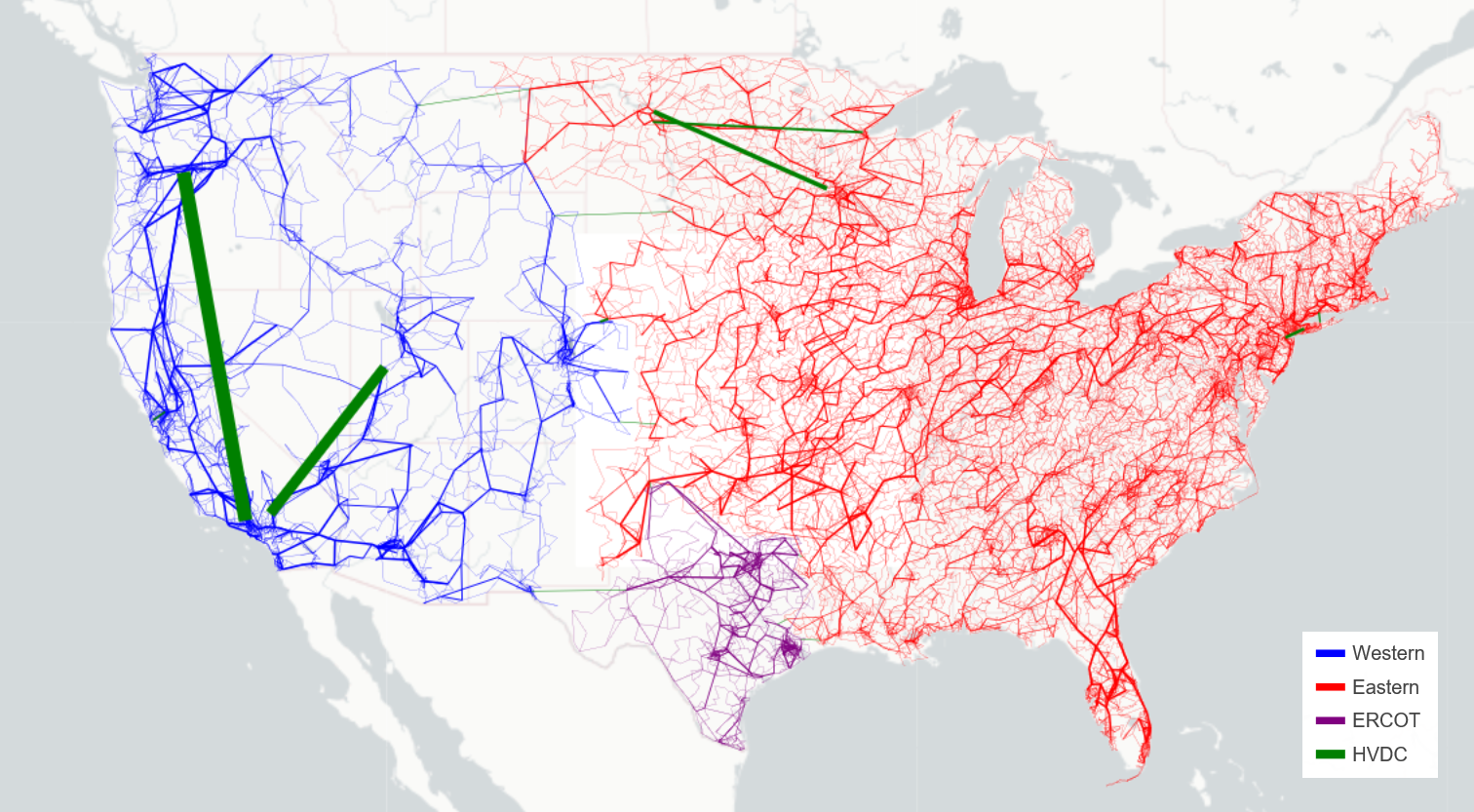}
    \caption{The high-resolution network topology used in the test system. Line thickness is proportional to power capacity.}
    \label{usmap}
\end{figure}

As a starting point for the network topology, the TAMU network was chosen for its high spatial resolution; it consists of 82,000 buses, 104,121 branches, and 13,419 generators in a fictitious configuration that is `realistic' and representative of true infrastructure but not contain any confidential information. This network is then updated with more recent generator capacities (described in Sections \ref{sec_fossil} to \ref{sec_hydro}), HVDC transmission lines (described in Section \ref{sec_hvdc}), and updated transmission capacity (described in Section \ref{sec_valid_network}). The final network is shown in Figure \ref{usmap}.

\subsection{Fossil Fuel Generators} \label{sec_fossil}

EIA Form 860 \cite{eia_860} was used to calculate the total operational generator capacity by fuel type for each state in the US (dividing states by interconnection when necessary). The capacity of the corresponding generators per state and fuel type in the TAMU network were scaled up or down to match the EIA capacities. When scaling generator capacity, the minimum power, ramp rates, and no-load costs were also scaled by the same value to maintain per-unit operational characteristics.

Within each state, the TAMU network features a range of cost curve coefficients for each fuel type (except for fuel oil generators, which are modeled as zero cost). However, the relative competitiveness of generators by state inadequately represents reality; for example, Wyoming coal generators were modeled as more expensive than many other states in the west despite their access to cheap Powder River Basin coal. To correct for this, average fuel prices for coal and natural gas by state and for fuel oil at the New York Harbor hub were obtained from historical EIA data \cite{eia_coal, eia_oil, eia_natgas} and paired with average generation heat rates by fuel \cite{eia_annual2018} to obtain average energy prices by fuel. Cost curves for all generators per fuel type in each state were scaled by a constant factor (or replaced in the case of fuel oil) such that capacity-weighted average electricity prices matched the historical data.

Emissions curves are obtained from heat rate curves by assuming constant rates of $\text{CO}_2$ per MMBTu of fuel burned \cite{eia_ghg}. These curves can be used to incorporate emissions into power system operations, for either an emissions-trading system or an emissions tax rate.

\subsection{Geothermal and Other Renewables}

Although geothermal generators are a relatively small share of generation capacity, they provide a more significant share of energy in California and Nevada due to their high capacity factors. These generators are not present in the TAMU network, so new generators are added to represent significant geothermal generators in these states. For California, capacity and generated energy per facility is available from the California Energy Commission \cite{cec_geo}. For Nevada, facility-level capacity data are available from NV Energy \cite{nv_energy}, state-level energy production data are available from the EIA \cite{eia_923}, and a constant capacity factor for all plants is assumed.

From inspection of generator dispatches from CAISO, geothermal generators run fairly constant generation rates throughout the year. Therefore, for each generator the maximum power rating is set to the average generation over the year and the minimum power rating is set to 95\% of the maximum power rating.
Geothermal facilities are assigned to network locations by correlating publicly available locations for the facilities with corresponding locations in the network with high voltage (100 kV or higher) buses. For `clusters' of geothermal generators (e.g. the Geysers, Salton Sea in California), the aggregate capacity is spread among several high voltage buses in the area, and the powerflow checked to ensure network feasibility and no undue local congestion.

Besides adding geothermal generators, new renewable generators are added when there were no utility-scale generators of that type in that zone (e.g. solar in the southeastern states). These generators are added at locations with good resource quality and adequate transmission capacity.

\subsection{Hydro Capacity} \label{sec_hydro}

The TAMU network includes capacities for Texas hydroelectric generators that do not adequately reflect current capacity. The capacities for this test system are adjusted using publicly available data from ERCOT to more closely align to real values \cite{ercot_cdr}. From conversations with TAMU, it was determined that some natural gas generators located near hydroelectric facilities were incorrectly coded as hydro in the TAMU network. After these generators were properly coded, capacities more closely matched data from EIA Form 860 \cite{eia_860}. Hydroelectric dams which are not present in the TAMU network are added at high voltage buses (100 kV or higher) near their true locations.

\subsection{HVDC Lines} \label{sec_hvdc}

Within each interconnection the TAMU network features only AC transmission lines, and a small number of HVDC lines are used to couple the interconnections. However, several significant HVDC lines in the US run entirely within interconnections. Therefore, additional intra-interconnection HVDC lines are added to the network by matching lat/long coordinates for true lines to coordinates of very-high voltage buses in the synthetic network. For example, the 3100 MW Pacific DC Intertie is modelled as running from a 765 kV bus in Oregon to a 345 kV bus in Southern California; although these voltages may not match the endpoint voltages of the true grid, these synthetic endpoints are connected to locations in the AC network with enough capacity to carry the full power that the DC line can provide.

\section{Hourly Time-Series Data}\label{sec_tsd}
\begin{figure}[h]
    \centering
    \includegraphics[width=\linewidth]{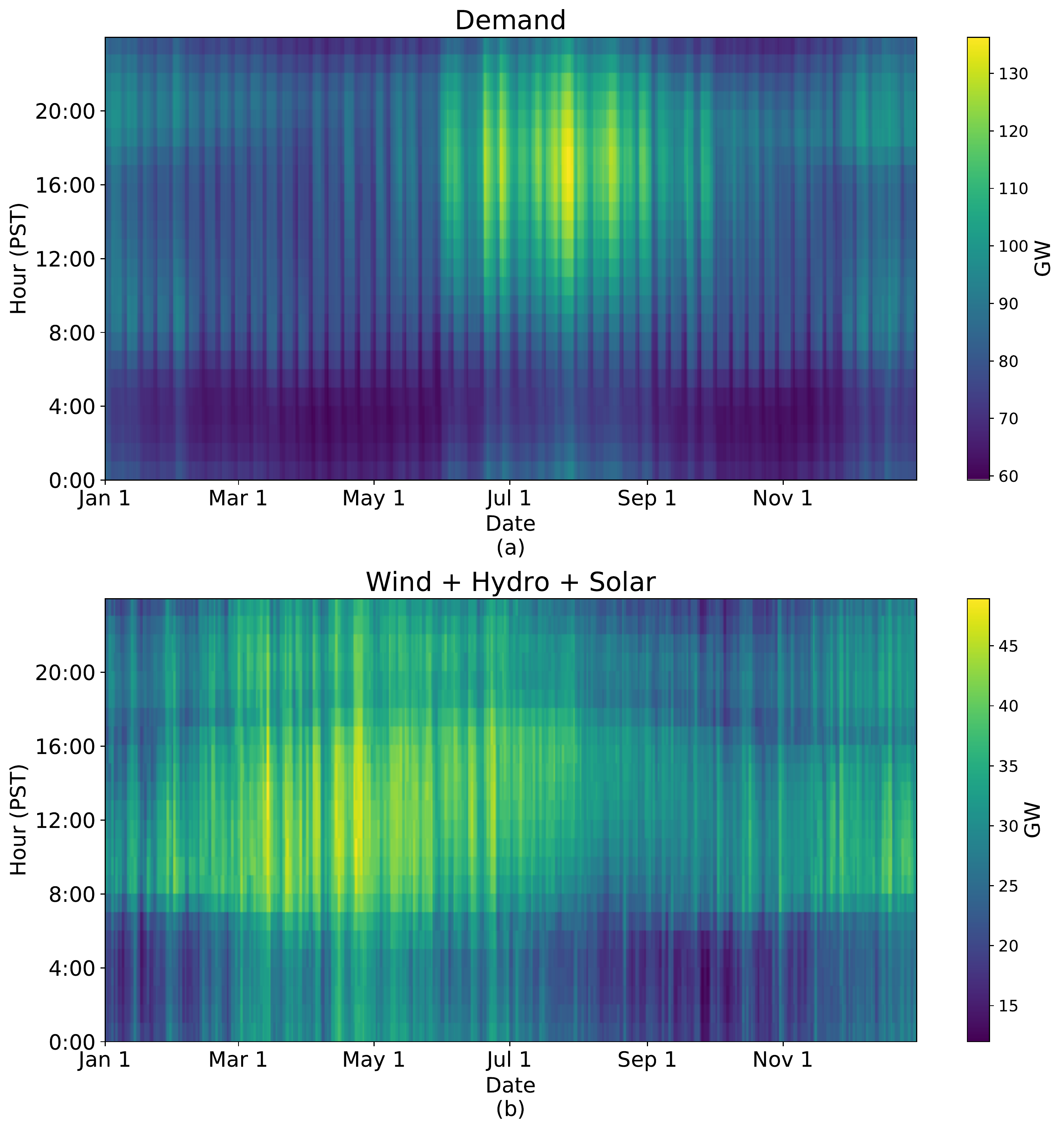}
    \caption{Hourly profiles of (a) demand and (b) solar, wind, and hydropower for the Western Interconnection.}
    \label{heatmap_both}
\end{figure}

%Hourly time series profiles representing each hour of the year 2016 are generated for demand, hydro, solar, and wind generation. The variability of the resulting demand and VRE profiles for the Western Interconnection are shown in Figures \ref{heatmap_demand} and \ref{heatmap_renewable} respectively.
Hourly time series profiles for 2016 are generated for demand buses and for hydro, solar, and wind generators. To avoid time zone and Daylight Saving Time discrepancies, all time series are standardized to UTC. Variability of the resulting demand and VRE profiles for the Western Interconnection is shown in Figure \ref{heatmap_both}. Interconnection demand peaks during late afternoons in summer while VRE generation peaks mid-day in spring, illustrating the time-alignment challenge of integrating large shares of VRE generation.

%\begin{figure}[h]
%    \centering
%    \includegraphics[width=\linewidth]{heatmap_demand}
%    \caption{heatmap demand}
%    \label{heatmap_demand}
%\end{figure}
%\begin{figure}[h]
%    \centering
%    \includegraphics[width=\linewidth]{heatmap_renewable}
%    \caption{heatmap renewable}
%    \label{heatmap_renewable}
%\end{figure}

\subsection{Demand Profiles}

Due to differing levels of data availability, different approaches are undertaken to construct demand profiles for buses in each interconnection.

\subsubsection{Western Interconnection}

The historical hourly load profile of each balancing authority area (BAA) in the U.S. is publicly available from the EIA \cite{eia_load}. Each county in the Western Interconnection is assigned to a BAA \cite{wecc_ba}, and total hourly load by BAA is distributed to each bus proportional to the population by ZIP code \cite{census}.

\subsubsection{ERCOT}

The historical hourly load profile of Texas weather zones is publicly available from ERCOT \cite{ercot_load}. The TAMU network has eight load zones which are geographically consistent with ERCOT's eight weather zones. Each load zone's demand profile is decomposed to bus-level load profiles proportional to the population of the ZIP code associated with each bus in the corresponding load zone.

\subsubsection{Eastern Interconnection}

The Eastern Interconnection is more complex than the Western Interconnection or ERCOT, and is represented by 70,000 nodes in the TAMU network compared to 10,000 for the Western Interconnection and 2,000 for Texas. As such, more effort was required to gather and clean demand data. Specifically:

\begin{enumerate}
    \item Hourly demand data were downloaded from the EIA \cite{eia_load}. Missing weekday data were replaced with an average of the corresponding hours in neighboring weekdays, while missing weekend data were duplicated from the remaining weekend day, or averaged from neighboring weekends when the entire weekend was unavailable.
    \item Anomalous demand data were identified as those with an hourly ramp rate magnitude of five standard deviations above average. These data points were discarded and linearly interpolated from non-anomalous hours.
    \item Counties were mapped to BAA via publicly available information or via direct contact with balancing authority representatives, and buses were mapped to counties using the Census Bureau's API \cite{census_api}.
    \item Finally, BAA hourly loads were decomposed to bus-level load profiles proportional to the population of the ZIP code associated with each bus in each BAA.
\end{enumerate}

\subsection{Hydro Profiles}

Compared to thermal generators, hydroelectric generators are significantly more constrained by hydrological cycles, and must be operated in accordance with agricultural, navigational, and ecological constraints. Monthly net generation is obtained from EIA Form 923 \cite{eia_923}, and historical hourly data from the largest twenty hydroelectric dams in the U.S. Army Corps of Engineers' Northwestern hydro system \cite{corp_army} are used to obtain a generic hydro profile to be applied to the generators in the Western Interconnection, with the exception of Wyoming and California. In Wyoming, the aforementioned method results in hourly generation that occasionally exceeds capacity, so monthly averages are used instead. In California, ISO data indicate that hydro profiles follow the net demand of of California \cite{caiso_outlook}, so this shape is applied to EIA energy totals.

Due to limited data availability for the Eastern Interconnection and Texas, hydro profiles are modeled as constant throughout the month in the Eastern Interconnection based on the average power from the EIA Form 923 data \cite{eia_923}. Future work will improve the fidelity of the Eastern hydro profiles.

\subsection{Solar Profiles}

Time-series solar irradiance for each solar plant in the TAMU network is obtained from the National Solar Radiation Database \cite{nrel_nsrdb}, which provides 1-hour resolution solar radiation data for the entire U.S. and a growing list of international locations on a 4~$\text{km}^{2}$ grid, spanning 1998-2006. Irradiance data are input into NREL's System Advisor Model (SAM) \cite{nrel_sam} to calculate power output. The PVWatts v5 model \cite{nrel_pvwatts} in SAM is used with all default values to obtain power profiles for each type of PV array (fixed open rack, 1-axis and 2-axis), and a weighted average of the different types is assigned to each solar plant in the TAMU network based on the overall ratio of installed capacity by interconnection. The prevalence of tracking systems is seen to be highly dependent on location, with the majority of capacity in the Eastern Interconnection utilizing fixed-tilt systems (only 33\% tracking) and the majority of capacity in the Western Interconnection and Texas utilizing single- or double-axis tracking systems (76\% for Western and 92\% for Texas). Finally, capacities of the solar plants in the TAMU network were scaled up by state such that the state totals matched EIA Form 860 data \cite{eia_860}.

\subsection{Wind Profiles}
Wind profiles are derived using the Rapid Refresh (RAP) \cite{noaa_rapid} numerical weather model. It covers North America and is comprised primarily of a numerical weather model and an analysis system to initialize that model. RAP provides latitudinal and longitudinal components (U and V) of the wind speed at 80~m above ground on a 13~$\text{km}^2$ resolution grid, for every hour ranging from May 2012 to date. Some data are missing in 2016 and hence are filled in using a simple procedure: extrema of the U and V components of the wind speed are calculated from all non-missing entries sharing the same location, same month, and same hour of the missing entry. Then, a U and V value are randomly generated within the respective derived ranges and used for imputation.

Wind speed is converted to power for all wind farms in the TAMU network using the IEC Class 2 power curve provided by NREL in the WIND Toolkit documentation \cite{nrel_wind}. Capacities of  wind farms in the TAMU network are scaled up by state such that state totals matched EIA Form 860 data \cite{eia_860}.

%\section{Simulation Platform}\label{sec_simplat}
%\input{04_Simulation_Platform}
\section{Validation Process}\label{sec_validation_process}
The high-resolution test system developed as described in Sections \ref{sec_static}-\ref{sec_tsd} was validated by running a multi-period optimal power flow using a DC powerflow approximation (MPDCOPF) with MATPOWER/MOST \cite{zimmerman2011, murillo2013}, with the multi-period dispatch chosen to capture the impact of ramp rate constraints. The 366 days of time-series data were broken into 61 individual simulations, each one representing 144 consecutive hours, and run consecutively. To bridge ramp-rate constraints across these boundaries, the optimal power output for each generator in the final hour of each simulation was used to constrain the active power of generators in the first hour of the next simulation. Any time the MPDCOPF was found to be infeasible, demand was reduced by 5\% for that simulation and the MPDCOPF for that time interval was re-launched with the reduced demand.

\subsection{Network Topology} \label{sec_valid_network}

When initially running the test system, unusual levels of congestion were observed in several branches in the Western Interconnection, and unusually severe load curtailment was observed in ERCOT. These observations prompted enhancements of the transmission network to better match the increased levels of demand and generation, compared to demand data that was originally developed based on census data from 2010 \cite{gegner2016}. Across all interconnections, transmission network capacity was increased until all demand can be feasibly served.

In the Western Interconnection, several branches labelled as \textit{Transformer} or \textit{TransformerWinding} exhibited especially high levels of congestion, with average shadow prices of congestion as high as \$119/MWh for one transformer in Idaho. High congestion shadow prices led to very low locational marginal prices (LMPs) in the surrounding upstream area, and therefore a large degree of renewable curtailment from solar and wind generators. Capacities of these congested branches were gradually increased in step sizes of around 100 MW until the LMPs at the renewable generators approached more reasonable averages. Overall, the capacities of only 9 branches were changed to reduce unrealistic congestion and curtailment.

In ERCOT, large increases in demand (primarily in the eastern half of the state) and large increases in wind generation (primarily in the western half of the state) were not reflected in the east-to-west transmission capacity, causing both load shedding in the east and wind curtailment in the west. To identify which branches were most critical to upgrade, a MPDCOPF problem with soft constraints on transmission capacity (with a large linear penalization factor in the objective) was used to identify which lines to upgrade, and by how much.

In all interconnections, many renewable generators are connected to the meshed grid via radial `spur' lines. When the capacity of these renewable generators are scaled up, the capacity of these spur lines is also increased to match the new renewable capacity. This ensures that unrealistic congestion does not exist on these spurs, as renewable generation investors will typically upgrade spur transmission capacity with renewable capacity, at least to the meshed grid.

\subsection{Cost Curve Revisions}

After running the model for a full year, generation totals by type and by state were compared against historical data from EIA Form 923 \cite{eia_923}. When generation quantities differed significantly from historical data, cost coefficients were adjusted on a state-by-state basis (by at most 5\%). In reality there are many mechanisms by which generation quantities are decided which are not captured by a transmission-constrained economic dispatch. However, modeling these mechanisms is outside the capability of most researchers; the cost-curve adjustments are a proxy which can yield similar results with less modeling and computational complexity.

\section{Validation Results}\label{sec:davail}
\begin{figure}[h]
    \centering
    \includegraphics[width=\linewidth]{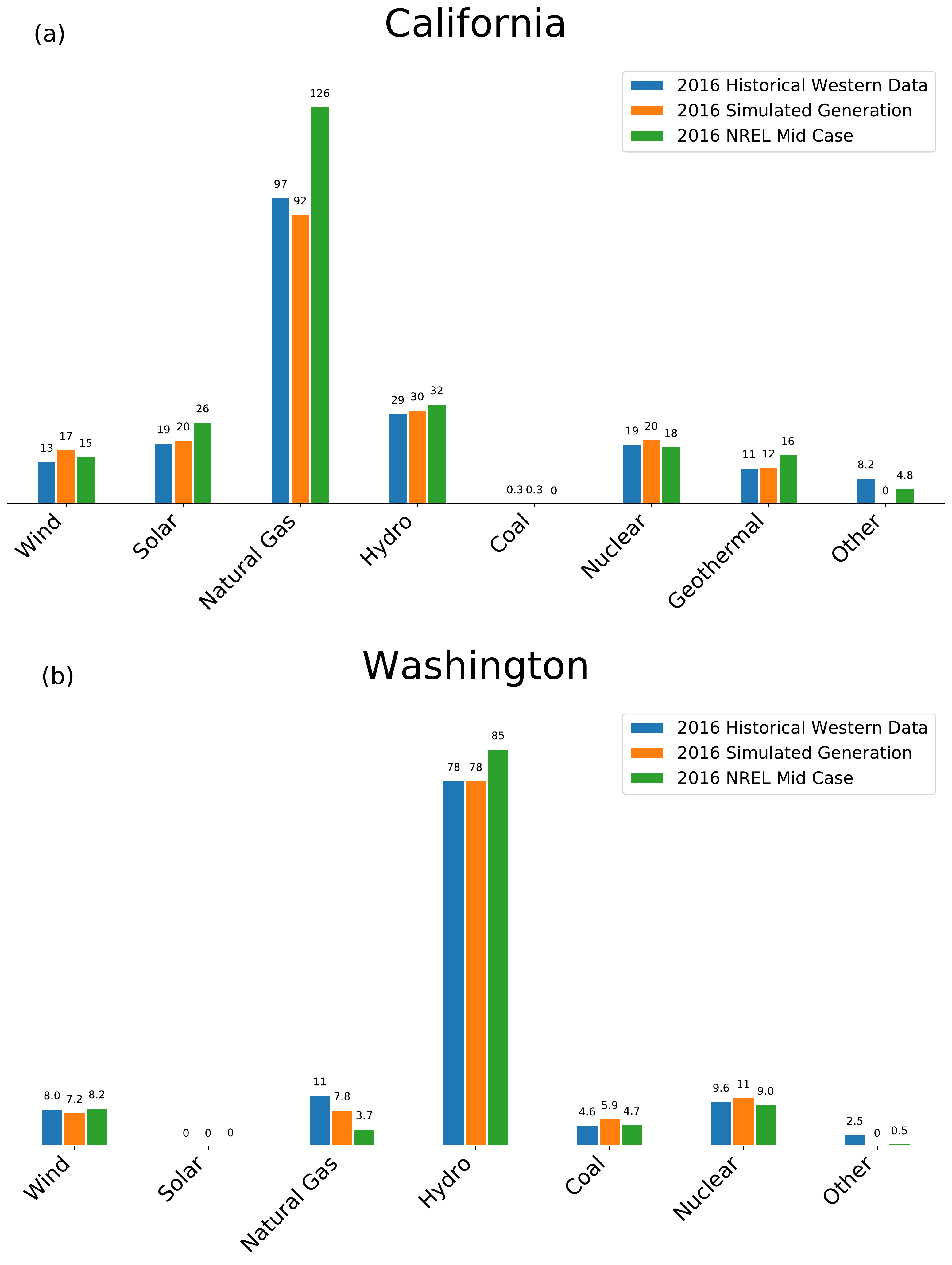}
    \caption{Total 2016 energy generation by fuel (TWh) in (a)~California and (b)~Washington.}
    \label{NREL_vs_Sim}
\end{figure}

To validate the results from the new dataset, the test system ran for an entire year as described in Section \ref{sec_validation_process} and results were compared to historical generation and to the NREL ReEDS Mid-Case Scenario results for 2016. Energy generation by fuel for California and Washington are shown in Figure \ref{NREL_vs_Sim}. Although the test system is still missing generators of less common types (e.g. biomass), for the included generator types the annual energy total matches fairly well with historical data, and compares favorably with the ReEDS results. Although ReEDS is more accurate for nuclear and wind generation, the new test system is more accurate for solar, natural gas, hydro, and geothermal generation; the overall energy error for the new system is smaller when measured by Euclidean distance or summed absolute errors. Total run-time is approximately 30 minutes for ERCOT, 3 hours for the Western Interconnection, and 20 hours for the Eastern Interconnection.

%\begin{figure}[h]
%    \centering
%    \includegraphics[width=\linewidth]{bars_ca}
%    \caption{Total energy consumption by fuel in California, 2016}
%    \label{bars_ca}
%\end{figure}

%\begin{figure}[h]
%    \centering
%    \includegraphics[width=\linewidth]{bars_wa}
%    \caption{Total energy consumption by fuel in Washington, 2016}
%    \label{bars_wa}
%\end{figure}

\section{Conclusions}\label{sec:conc}
A set of methods has been presented to develop high-resolution, publicly available test systems with coverage of continental-scale power systems. Future work will include development of input datasets under various technical and political scenarios which will impact projected demand profiles, transmission capacities, and fuel prices.

The primary goals of developing enhanced test systems as described herein are \textit{a)} to support more accurate modeling of current U.S. infrastructure, and \textit{b)} to support development and testing of high-resolution capacity expansion models (CEMs). Many commercial CEMs and datasets exist, and several open-source or open-access tools are in active development (e.g. ReEDS \cite{nrel_reeds}, SWITCH 2.0 \cite{johnston2019}, and PyPSA \cite{pypsa}), but currently there appears to be no fully-open source CEM tool with associated high-resolution, non-confidential data for the U.S.

The authors are currently developing an open-source CEM which can make use of high-resolution data to determine the optimal combination of equipment (generation, transmission, energy storage, distributed energy resources, etc.) for deep decarbonization efforts. Design goals include a full set of research-grade features and a flexible user interface accessible to users with a wide range of technical backgrounds. 

\bibliographystyle{IEEEtran}
\bibliography{bibliography}

\end{document}